\documentclass[12pt]{article}
\usepackage{amssymb,amsmath,latexsym}

\def\proof{{\it Proof}.\ }

\def\wbull{\hfill\vrule height .9ex width .8ex depth -.1ex}

\newtheorem{formula}{}[section]
\newtheorem{proposition}[formula]{Proposition}
\newtheorem{definition}[formula]{Definition}
\newtheorem{corollary}[formula]{Corollary}

\newtheorem{lemma}[formula]{Lemma}
\newtheorem{theorem}[formula]{Theorem}

\def\thrm{\begin{theorem}}
\def\thrml#1{\begin{theorem}\label{#1}}
\def\ethrm{\end{theorem}}
\def\prp{\begin{proposition}}
\def\prp#1{\begin{proposition}\label{#1}}
\def\eprp{\end{proposition}}
\def\dfntn{\begin{definition}}
\def\dfntnl#1{\begin{definition}\label{#1}}
\def\edfntn{\end{definition}}
\def\nmrt{\begin{enumerate}}
\def\enmrt{\end{enumerate}}

\def\qtn{\begin{equation}}
\def\qtnl#1{\begin{equation}\label{#1}}
\def\eqtn{\end{equation}}
\def\lmm{\begin{lemma}}
\def\lmml#1{\begin{lemma}\label{#1}}
\def\elmm{\end{lemma}}
\def\crllr{\begin{corollary}}
\def\crllrl#1{\begin{corollary}\label{#1}}
\def\ecrllr{\end{corollary}}
\def\css{\begin{cases}}
\def\ecss{\end{cases}}

\title{ \bf{There exists no distance-regular graph with intersection array 
$\{55,36,11;1,4,45\}$}}
\author{
Alexander Gavrilyuk
\thanks{
Partially supported by the Russian Foundation for Basic Research (project no. 08-01-00009).}\\[-1pt]
\small Ural Division of the Russian Academy of Sciences\\[-3pt]
\small Institute of Mathematics and Mechanics\\[-3pt]
\small ul.~S.~Kovalevskoi~16, Yekaterinburg, 620990 Russia\\[-3pt]
{\tt \small alexander.gavriliouk@gmail.com}\\[-3pt]}

\begin{document}

\maketitle

\begin{abstract}
We prove that a distance-regular graph with intersection array 
$\{55,36,11;1,4,45\}$ does not exist. This intersection array
is from the table of feasible parameters for distance-regular graphs
in "Distance-regular graphs"\ by A.E. Brouwer, A.M. Cohen, A. Neumaier.
\medskip

\end{abstract}
\newpage

\setcounter{section}{1}


The purpose of this short paper is to show that a distance-regular 
graph with intersection array $\{55,36,11;1,4,45\}$ does not exist. 
We use the notations and definitions given in~\cite{BCN}.
In particular, for a graph $\Gamma$ and a vertex $x\in \Gamma$, 
we denote the neighborhood of $x$ by $\Gamma(x)$.

Let us first recall some needed results.
Let $\Gamma$ be a distance-regular graph with intersection array 
$\{b_0,b_1,\ldots,b_{d-1};1,c_2,\ldots,c_{d}\}$ and 
$x$ be an arbitrary vertex of $\Gamma$.
The neighborhood of $x$ is a regular graph with degree $a_1:=b_0-b_1-c_1$.
The adjacency matrix of $\Gamma$ has exactly $d+1$ distinct eigenvalues, 
say, $b_0=\theta_0>\theta_1>\ldots>\theta_d$.

\lmml{KPineq}
If the graph $\Gamma(x)$ contains a coclique of size $c\ge 2$, then
$$c_2-1\ge \displaystyle {{\frac {c (a_1+1)-b_0}{{c\choose 2}}}}.$$
\elmm
\proof See~\cite{KP}.\wbull

\lmml{Tineq}
The neighborhood of $x$ is a graph with the second largest eigenvalue 
$\le \displaystyle { -{\frac {b_1}{\theta_d+1}}-1}$ 
(here the second largest eigenvalue is taken to be 
the degree $a_1$ in case of disconnected $\Gamma(x)$).
\elmm
\proof See~\cite[Th. 4.4.3]{BCN}.\wbull

\lmml{Brineq}
Let $\Delta$ be a connected graph on $v$ vertices. 
If $\Delta$ is regular with degree $k<v-1$ and is not an odd cycle, 
then $\Delta$ contains a coclique of size 
$\ge v/k$.
\elmm
\proof It follows from~\cite[Prop. 4.2.1]{BH} that the chromatic number 
$\chi(\Delta)\le 1+k$ with equality if and only if $\Delta$ is complete 
or an odd cycle. The size of the largest coclique in $\Delta$ is 
at least 
$v/\chi(\Delta)$.\wbull
\medskip

\thrml{T}
The array $\{55,36,11;1,4,45\}$ cannot be realized as 
the intersection array of a distance-regular graph.
\ethrm
\proof Let $\Gamma$ be a distance-regular graph with intersection array 
$\{55,36,11;1,4,\\45\}$ (see~\cite[p. 429]{BCN}). Then the adjacency matrix 
of $\Gamma$ has exactly $4$ distinct eigenvalues: $55,19,-1$ and $-5$.

Let $x$ be a vertex of $\Gamma$. By Lemma \ref{Tineq}, 
the second largest eigenvalue of $\Gamma(x)$ is at most 
$36/4-1<a_1=18$. Hence, 
the graph $\Gamma(x)$ is connected. By Lemma \ref{Brineq}, 
the graph $\Gamma(x)$ contains a coclique of size 
$\ge 55/18>3$.
Now, by Lemma \ref{KPineq}, we have 
$4-1\ge \displaystyle {{\frac {4(18+1)-55}{{4\choose 2}}}}=3.5$, 
which is impossible. The theorem is proved.
\bigskip


\end{document}